\pdfoutput=1
\RequirePackage{ifpdf}
\ifpdf 
\documentclass[pdftex]{sigma}
\else
\documentclass{sigma}
\fi

\begin{document}
\allowdisplaybreaks

\newcommand{\arXivNumber}{1907.05610}

\renewcommand{\thefootnote}{}

\renewcommand{\PaperNumber}{070}

\FirstPageHeading

\ShortArticleName{Holomorphic Distributions and Connectivity by Integral Curves of Distributions}

\ArticleName{Holomorphic Distributions and Connectivity\\ by Integral Curves of Distributions\footnote{This paper is a~contribution to the Special Issue on Algebra, Topology, and Dynamics in Interaction in honor of Dmitry Fuchs. The full collection is available at \href{https://www.emis.de/journals/SIGMA/Fuchs.html}{https://www.emis.de/journals/SIGMA/Fuchs.html}}}

\Author{Vladimir A.~ZORICH}

\AuthorNameForHeading{V.A.~Zorich}

\Address{Division of Mathematical Analysis, Faculty of Mechanics and Mathematics,\\
Lomonosov Moscow State University, GSP-1, Leninskie Gory, Moscow, 119991, Russia}
\Email{\href{mailto:vzor@mccme.ru}{vzor@mccme.ru}}

\ArticleDates{Received July 19, 2019, in final form September 14, 2019; Published online September 19, 2019}

\Abstract{It is known that the classical Frobenius theorem on conditions of integrability for distributions of planes can be extended to the case of complex holomorphic distributions. We show that an alternative criterion for integrability, namely, non-connectivity, discovered (or at least, marked and explicitly formulated) by Carath\'eodory in relation to classical thermodynamics, also admits a holomorphic formulation.}


\Classification{53D10; 32B15; 80A05}

\renewcommand{\thefootnote}{\arabic{footnote}}
\setcounter{footnote}{0}

\section{Introduction}

We consider a holomorphic version of a theorem of Carath\'eodory initiated by his studies related to mathematical aspects of classical thermodynamics. The generalization which we are going to discuss, like the Carath\'eodory theorem itself, is easily provable by modern tools of analysis and geometry. However, an interesting, useful and important feature is rather a remarkable duality of connectivity of the space of states of a thermodynamic system, or in any abstract space, by admissible (tangent, contact) paths for a distribution of hyper\-pla\-nes in the space, and integrability of this dis\-tri\-bu\-tion itself.

In more detail, in classical equilibrium thermodynamics one often has to consider adiabatic transitions of a thermodynamic system from one thermodynamic state to another. The adiabatic condition in mathematical language means the requirement that the corresponding transition path goes along the kernels (zero spaces) of a special 1-form, i.e., it is constrained to be tangent to a certain distribution of hyperplanes. Justifying (mathematically formalizing) the so-called second law of thermodynamics, generating the entropy function of a thermo\-dy\-na\-mic state, Carath\'eodory \cite{Caratheodory,Caratheodory+} investigated the possibility of adiabatic transitions between thermodynamic states. In mathematical terms he showed that the absence of local connectivity of points in space by paths tangent to the distribution of hyperplanes is equivalent to integrability of the distribution. Thus, the question of connectivity (non-connectivity) turned out to be related with the issue of integrability of the distribution given by a 1-form (and if we omit some details, it turned out to be equivalent to it). But the problem of integrability of a distribution was already solved by the well-known Frobenius theorem. Thus, the problem of connectivity reduces to the already solved problem.

Note that the first significant observations in this direction were already made by Poin\-ca\-r\'e~\cite{Poincare}, also in connection with thermodynamics.

The reference \cite{Born} is a nice extended exposition by Born of Carath\'eodory's work, possibly more accessible to physicists.
A mathematical generalization of Carath\'eodory's theorem for distributions of arbitrary codimension (not only to distributions of hyperplanes) is now known as the Rashevskiy--Chow theo\-rem~\cite{Chow, Rashevskiy}.

\section{Integrability of distributions}

Consider a smooth field (distribution) of hyperplanes in the space ${\mathbb R}^n$. Such a distribution can have integral surfaces, which are surfaces whose tangent planes at each point coincide with the corresponding planes of the distribution. If the space is locally foliated by such surfaces, then such a distribution of hyperplanes is called an integrable distribution.

Generic vector fields are always integrable (this is one of fundamental initial facts of the theory of ordinary differential equations). For distributions of higher dimension this is not always the case. For hyperplane distribution defined by kernels (zeros) of a non-degenerate 1-form $\omega$ a necessary and sufficient condition of local integrability of the distribution is the condition $\omega \wedge {\rm d}\omega \equiv 0$. This is the Frobenius theorem~\cite{Cartan}, or rather one of its versions.

Below we will be interested in holomorphic distributions of complex hyperplanes in the complex space. It can be verified that the local integrability condition expressed by the Frobenius theorem is also applicable when the distribution is given by kernels (zeros) of a holomorphic 1-form $\omega$. (See, for example,~\cite{Nierenberg} or proofs of the classical Frobenius theorem in~\cite{Cartan} or~\cite{Arnold}.)

\section{Connectivity by paths tangent to distributions}

Given a plane distri\-bu\-tion in space, we are allowed to move from one to another point of the space subject to the constraint that at each moment the velocity vector of motion has to lie in the corresponding plane of the distribution. In this case the path is said to be tangent (or subordinated) to the distribution, or that it is an admissible path or an integral curve for the distribution.

In mechanics this is a typical example of motion in the presence of nonholonomic constraints. In thermodynamics it may be a problem of the adiabatic transfer of a thermodynamic system from one equilibrium state to another. In control theory one looks for (im)possibility of reaching one point of space from another by means of a controlled admissible path for a given distribution, determined by the problem itself.

If the distribution is integrable, then obviously no admissible path can leave the integral surface of the distribution to which the starting point of the path belongs. No point lying outside of the integral surface can be connected by an admissible path with any point of this integral surface.

So, if there is integrability of a distribution of hyperplanes, then, certainly, there is no local connectivity of points of space by admissible paths. (For instance, in a thermodynamic system there cannot be an adiabatic transition through equilibrium states between two equilibrium states having different entropy values.)

A nontrivial observation of Carath\'eodory, initiated by his work on the ma\-the\-ma\-ti\-cal formalization of classical thermodynamics, is the above-mentioned statement that
in the case of a hyperplane distribution the inverse theorem is valid: if there is no local connectivity, then there is local integrability.

In short, although not quite accurately, there is an alternative: either there is local connectivity and there is no local integrability, or there is no local connectivity, and then there is local integrability.

Since the question of the local integrability of the distribution of hyper\-pla\-nes is solved by the Frobenius theorem, Carath\'eodory's theorem provides an effective tool to verify the possibility of connecting points (states of the system) by means of paths admissible for a given distribution.

\section{Main statement}

Now we can turn to the main subject of this note. We are going to show that not only the holomorphic analogue of the Frobenius theorem is valid, which was already mentioned above, but the holomorphic analogue of the Carath\'eodory theorem takes place as well.

We will consider the simplest nontrivial case, the distribution of complex planes in the complex space~${\mathbb C}^3$.

Let our distribution be given as zeros of the standard 1-form $ x{\rm d}y + {\rm d}z $. Instead of the form $ x{\rm d}y + {\rm d}z $ at the moment it is technically more convenient for us to consider the linearly equivalent form $ x{\rm d}y - y{\rm d}x - {\rm d}z $.

\looseness=-1 Let $ x = t$, $y = t^2 $, and $ z = \frac{1}{3} t^3 $. Since ${\rm d}z = x{\rm d}y - y{\rm d}x $, we have an admissible holomorphic curve parametrized by the complex parameter $ t $. The function $ z = \frac{1}{3} t^3 $ for some value of $t = t_1 $ can take any prescribed value~$ z_1 $. At the complex moment~$ t_1 $ the curve $ \big(x = t, y = t^2$, $z = \frac{1}{3} t^3\big) $ passes through the point $ (x_1, y_1, z_1) = \big(t_1, t^2_1,\frac{1}{3} t^3_1\big) $. By an additional transition along the admissible complex line $ (x = x_1 t, y = y_1 t, z = z_1) $ one can move to the point $(0, 0, z_1)$. But once we managed to get from the point $(0, 0, 0)$ to any point of the form $ (0, 0, z_1) $ by an admissible holomorphic curve, then, using the same construction, we can reach any point of the space ${\mathbb C}^3$.

Thus, in principle, it is possible to pass from any point of the space $ {\mathbb C}^3 $ to any other point along holomorphic curves subjected to the distribution defined in $ {\mathbb C}^3 $ as $ \{\ker(x{\rm d}y - y{\rm d}x - {\rm d}z) \}$ or as $\{\ker (x{\rm d}y + {\rm d}z) \}$.

In the above discussion a holomorphic broken line has no more than three links. One can carry out this reasoning more carefully, and instead of a piece-wise holomorphic curve use a single holomorphic curve connecting a selected pair of points. For instance, in order to connect in $ {\mathbb C}^3 $ the origin $(0, 0, 0)$ with an arbitrary point $ (x_1, y_1, z_1) $ by a holomorphic curve admissible for the distribution $ \{\ker(x{\rm d}y - y{\rm d}x - {\rm d}z) \}$, one can consider the curve $ \big(x = x_1 t, y = y_1 t^2 - ct(t-1), z = \frac 13 x_1 y_1 t^3 - \frac 13 c x_1 t^3\big) $ by choosing the constant $ c $ satisfying the condition $ z_1 = \frac 13 x_1 y_1 - \frac 13 c x_1 $ (for $ t = 1 $). This linear equation is solvable if $ x_1 \neq 0 $. If $ x_1 = 0 $ and $ y_1 \neq 0 $, then we can do this procedure by interchanging the variables. If $ x_1 = y_1 = 0 $ the task becomes extremely concrete: starting from the origin we have to reach the point $ (0,0, z_1) $ by holomorphic curve admissible for the given distribution.

Note that the quasi-homogeneous dilation $ (x, y, z) \mapsto \big(\alpha x, \alpha y, \alpha^2 z\big) $ pre\-ser\-ves our distribution, therefore, it suffices to check local connectivity in the neigh\-bour\-hood of the origin.

For the general holomorphic distribution of complex codimension 1 the Darboux theorem holds: by the suitable holomorphic change of local co\-or\-di\-na\-tes the form generating the distribution can be reduced to the normal form $ x{\rm d}y + {\rm d}z $ (see~\cite{Arnold}, for the proof of the classical Darboux theorem, and~\cite{Alarc}, for the proof in the case of a holomorphic distribution). Thus, ``generic'' holomorphic distributions are completely non-integrable in the sense that locally they allow transitions between points of space along holomorphic curves admissible for the distribution.

In the case of the space $ {\mathbb C}^3 $, in contrast to the general multidimensional case, holomorphic curves integral (admissible) to the form $ x{\rm d}y + {\rm d}z $ are holomorphic Legendrian integral manifolds (of maximal possible dimension). But all such varieties (in~$ {{\mathbb C}^{2n+1}} $ or in~$ {{\mathbb R}^{2n+1}} $) are described explicitly as follows
\begin{gather*}
y_I = \frac{\partial S}{\partial x_I}, \qquad x_J = - \frac{\partial S}{\partial y_J}, \qquad z = S - x_I \frac{\partial S}{\partial x_I},\end{gather*}
by means of the generating function $ S $ (see \cite[Appendix 4]{Arnold} on contact structures). Here $ S = S (x_I, y_J) $, and $ I \cup J $ is any partition of the indices $ (1,\dots, n) $ when we consider the space $ {\mathbb R}^{2n + 1} $ or the space $ {\mathbb C}^{2n + 1} $. In the case of $ {\mathbb C}^3 $ the function $ S $ depends only on one of the variables $ x, y $. This allows one (by selecting the function $ S $ as a polynomial $ ax^2 + bx^3 $ or as a polynomial $ ay^2 + by^3 $) directly verify that any point in the neighborhood of the origin can be connected to the origin by a holomorphic curve integral for the distribution $ \{ \ker(x{\rm d}y + {\rm d}z)\} $.

Now we will try to obtain an inverse theorem, that of the Carath\'eo\-do\-ry theorem type: if a holomorphic distribution of complex codimension~1 does not allow connections between any points in any neighborhood of the space by means of holomorphic curves, admissible for the distribution, then the distribution is integrable.

So, we will assume that in any neighborhood of every point $ p $ there are points inaccessible from $ p $ by admissible holomorphic curves. We fix some starting point $ p $ together with its certain neighborhood. The idea of the further heuristic reasoning could be as follows.\footnote{This idea can be seen already in Poincar\'e \cite[Section 193e]{Poincare}. Starting with the adiabatic inaccessibility of some equilibrium states of the thermodynamic system, Poincar\'e (and then Carath\'eodory) derives integrability of the corresponding differential 1-form of heat influx and obtains one of the most important characteristics of an equilibrium state of a thermodynamic system, the entropy. Levels of the entropy function (adiabats, isoentropes) integrate the distribution given by the heat influx 1-form.}

The set of points of this neighborhood inaccessible (unreachable) from the point $ p $ is an open set (since if accessible points converge to a point, then the limiting point is also accessible; of course, it requires further mathematical justification, although from a physical point of view it is almost a tautology).

Consider the boundary $ \Gamma $ of the set of inaccessible points. Our starting point $ p $ (the center of the neighborhood where events occur) must lie on~$ \Gamma $. The boundary $ \Gamma $ is expected to be regular, due to the smoothness of the distribution (for example, it cannot contain conic points and corners).

The specific plane of the distribution attached at the point~$ p $ cannot be transversal to~$ \Gamma $, as otherwise it would be possible to enter the region of inaccessible points along this plane. Therefore, the considered plane of the distribution turns out to be tangent to~$\Gamma $.

We would like to show that in the case of the holomorphic distribution we started with, $ \Gamma $ has complex codimension 1, and therefore $ \Gamma $ turns out to be an integral surface of the initial holomorphic distribution.

Note that the holomorphic property of the distribution must be used substantially. Indeed, consider, for example, the foliation of the space $ {{\mathbb C}^n} $ by concentric spheres and look at the distribution of complex tangents to the spheres. This distribution is not holomorphic (it is parameterized by spheres). If it were holomorphic, then in the case of $ {{\mathbb C}^2} $ it would even be holomorphically straightened and, like a vector field, it would be integrable. Within each of the spheres, the distribution of their complex tangents is completely non-integrable: there is possible to connect any two points of the sphere along real one-dimensional paths admissible for the distribution. However, it is impossible to leave the sphere along such a path. Therefore in this setting in the neighborhood of each point of the whole space there are points that cannot be reached even along real one-dimensional paths admissible for the distribution.

Now we return to a holomorphic distribution in~$ {\mathbb C}^3 $. A somewhat more detailed study of the boundary set~$ \Gamma $ leads to the desired conclusion that $ \Gamma $ is a complex hypersurface and it is integral for the initial holomorphic distribution of hyperplanes.
Instead, we will replace such a description of the set $ \Gamma $ by the following shorter argument\footnote{It was proposed by F.~Forstneri\v{c} and is published with his kind permission.} using that distribution is holomorphic, while the space is three-dimensional.
In the space $ {\mathbb C}^3 $ the 3-form $ \omega \wedge {\rm d}\omega $ is determined by only one functional coefficient. If the coefficient is identically equal to zero, then $ \omega \wedge {\rm d}\omega \equiv 0 $ and the Frobenius condition for integrability of the distribution is fulfilled.

On the other hand, the non-connectivity assumption implies that there are no points where this coefficient differs from zero. Indeed, if $ \omega \wedge {\rm d}\omega \not \equiv 0 $, then the equation \makebox{$ \omega \wedge {\rm d}\omega = 0 $} determines a~complex hypersurface in $ {\mathbb C}^3 $ outside of which the holomorphic 2-form $ {\rm d}\omega $ must be non-degenerate on~$ \ker\omega $.

But in a neighborhood of points where the form $ {\rm d}\omega $ is nondegenerate, the form $ \omega $ can be written in the form $ \omega = x{\rm d}y + {\rm d}z $, thanks to the Darboux theorem. Hence, as we saw above, there appears accessibility of all points of that neighborhood by admissible paths, while this contradicts the assumption. Thus, the following theorem holds.
\begin{theorem*}Consider a holomorphic distribution of complex hyperplanes in the space $ {\mathbb C}^3 $. If in any neighborhood of each point of the space there are points inaccessible from the center of the neighborhood by holomorphic curves admis\-si\-ble for the distribution, then the distribution is integrable.
\end{theorem*}
Considering holomorphic distributions of complex hyperplanes we confine ourselves to the first nontrivial case.

To complete the general case consider the form $ \omega \wedge ({\rm d} \omega)^k $, where $k$ is minimal of the natural numbers such that $ \omega \wedge ({\rm d} \omega)^k \equiv 0 $. If $k =1$, we have the integrability condition. On the other hand, $k$ must be equal to~1, otherwise the non-connectivity condition will be violated. It follows from the Darboux theorem applied to the form $ \omega \wedge ({\rm d} \omega)^{k-1} $ in the neighbourhood of the point where the form is not equal to zero.
(The following example can explain the above reasoning. If in the space $ {\mathbb C}^5 = {\mathbb C}^3 \times {\mathbb C}^2 $ instead of the form $ x_1 {\rm d}y_1 + x_2 {\rm d}y_2 + {\rm d}z $ one considers the form $ \omega = x_1 {\rm d}y_1 + {\rm d}z $, then $ \omega \wedge ({\rm d} \omega)^2 \equiv 0$ but $ \omega \wedge {\rm d} \omega \neq 0 $.)

Note that since the theorem is actually local, it suggests the similar considerations for an arbitrary complex manifold.

\section{Final comments}

The idea of Poincar\'e mentioned above (very natural in the context of thermodynamics) is particularly attractive, since it directly and explicitly relates the problem of connectivity of the points of the space by admissible paths and (non)integrability of the corresponding distri\-bu\-tion, presenting the integral surface $\Gamma$.

However, the detailed proof of the fact that the manifold $\Gamma$ is a holomorphic hypersurface turns out to be rather intricate. Thus, it is worth saying a few words about another approach to the relation between integrability of the distribution and connectivity of points of the space by integral paths of the distribution. We mean the following proof (essentially being an illustration to the Frobenius theorem), which works both in real and complex cases and for different dimensions of the space under consideration.

Above we have already mentioned Frobenius criterion $ \omega \wedge {\rm d}\omega \equiv 0 $ for the distribution $\{\ker\omega\}$ to be an integrable one. In turn, this condition is equivalent to the fact that the differential $ {\rm d}\omega $ of the 1-form $ \omega $ vanishes on any pair of vector fields tangent to the planes of that distribution.

But according to Cartan's formula, known in calculus of differential forms (see, e.g.,~\cite{Cartan}), if $ X $ and $ Y $ are any smooth vector fields and $ \omega $ is a 1-form on a manifold, then ${\rm d}\omega(X, Y) = X \omega(Y) - Y \omega(X) - \omega([X, Y])$, where $ X \omega(Y) $ and $ Y \omega(X) $ are Lie derivatives of the functions $ \omega(Y) $ and $ \omega(X)$ along the fields~$ X $ and~$ Y $, respectively, and $ [X, Y] $ is the Lie bracket (commutator) of these vector fields.

If the vector fields belong to the planes of the distribution, then $ \omega(Y) \equiv 0 $ and $ \omega(X) \equiv 0 $, and the equality $ {\rm d}\omega(X, Y) = 0 $ depends on whether the commutator $ [X, Y] $ belongs to the planes of the distribution.

If the commutator is not tangent to the distribution, then, as is well known, by going along such fields one can move in directions transversal to the planes of distribution
(see~\cite{Arnold,Gromov,Montgomery}).

Thus, Cartan's formula provides a direct connection between inte\-gra\-bi\-li\-ty of a hyperplane distribution in space and connectivity of points of the space by integral paths of the distribution. Namely, in the case of a typical distribution the brackets of vector fields spanning the distribution of hyperplanes can be transversal to the hyperplanes of distribution. This provides local connectivity of points of the space by paths admissible for the distribution.

For distributions of higher codimension, and even in degenerate cases of hyperplane distributions, instead of integral surfaces of distribution, there could appear such foliations of the space that the distribution is completely non-integrable and the connectivity of points by admissible paths is possible within each leaf of the foliation, but the transfer between the leaves is impossible. We have already seen such a phenomenon in examples above. Singular cases, of course, require special consideration and description. They were not presented in the initial papers by Poincar\'e and Carath\'eodory mentioned above and cited in the references.

In the introduction we briefly explained a relation of the issue discussed above, as well as of Carath\'eodory's theorem, with classical thermo\-dy\-na\-mics. One can read more about that relation, e.g., in~\cite{Zorich1} and~\cite{Zorich2}.

\vspace{-1.5mm}

\subsection*{Acknowledgements}

I am very grateful to S.Yu.~Nemirovski for fruitful discussions and for the reference~\cite{Alarc}, where the authors not only report their impressive results, but provide full proofs of several facts of mathematical folklore (e.g., holomorphic versions of the Darboux theorem), which they use, and which we also needed above. I am also indebted to F.~Forstneri\v{c} for his argument used above, and for pointing out the new reference~\cite{AlarcForst}. Special thanks to B.~Khesin who carefully read my source text and corrected not only my English but my math as well. Additional thanks go to the referees for valuable remarks and suggestions.



\vspace{-3mm}

\pdfbookmark[1]{References}{ref}
\LastPageEnding

\end{document}